
\magnification = \magstep0
\overfullrule = 2pt
\vsize = 525dd
\hsize = 27cc 
\topskip = 13dd
\hoffset = 1.8cm
\voffset = 1.8cm


\font\bfbig = cmbx10 scaled \magstep2   
\font\partf = cmbx10 scaled \magstep1	

\font\eightrm = cmr10 scaled 800        
\font\sixrm = cmr7 scaled 850
\font\fiverm = cmr5
\font\eighti = cmmi10 scaled 800
\font\sixi = cmmi7 scaled 850
\font\fivei = cmmi5
\font\eightit = cmti10 scaled 800
\font\eightsy = cmsy10 scaled 800
\font\sixsy = cmsy7 scaled 850
\font\fivesy = cmsy5
\font\eightsl = cmsl10 scaled 800
\font\eighttt = cmtt10 scaled 800
\font\eightbf = cmbx10 scaled 800
\font\sixbf = cmbx7 scaled 850
\font\fivebf = cmbx5

\font\fivefk = eufm5                    
\font\sixfk = eufm7 scaled 850
\font\sevenfk = eufm7
\font\eightfk = eufm10 scaled 800
\font\tenfk = eufm10

\newfam\fkfam
	\textfont\fkfam=\tenfk \scriptfont\fkfam=\sevenfk
		\scriptscriptfont\fkfam=\fivefk

\def\eightpoint{%
	\textfont0=\eightrm \scriptfont0=\sixrm \scriptscriptfont0=\fiverm
		\def\rm{\fam0\eightrm}%
	\textfont1=\eighti  \scriptfont1=\sixi  \scriptscriptfont1=\fivei
		\def\oldstyle{\fam1\eighti}%
	\textfont2=\eightsy \scriptfont2=\sixsy \scriptscriptfont2=\fivesy
	\textfont\itfam=\eightit \def\it{\fam\itfam\eightit}%
	\textfont\slfam=\eightsl \def\sl{\fam\slfam\eightsl}%
	\textfont\ttfam=\eighttt \def\tt{\fam\ttfam\eighttt}%
	\textfont\bffam=\eightbf \scriptfont\bffam=\sixbf
		\scriptscriptfont\bffam=\fivebf \def\bf{\fam\bffam\eightbf}%
	\textfont\fkfam=\eightfk \scriptfont\fkfam=\sixfk
		\scriptscriptfont\fkfam=\fivefk
	\rm}
\skewchar\eighti='177\skewchar\sixi='177\skewchar\eightsy='60\skewchar\sixsy='60

\def\fk{\fam\fkfam}
\def\petit{\eightpoint}%


\def\today{\ifcase\month\or
	January\or February\or March\or April\or May\or June\or
	July\or August\or September\or October\or November\or December\fi
	\space\number\day, \number\year}
\def\newline{\hfil\break}
\def\framebox#1{\vbox{\hrule\hbox{\vrule\hskip1pt
	\vbox{\vskip1pt\relax#1\vskip1pt}\hskip1pt\vrule}\hrule}}


\newtoks\RunAuthor\RunAuthor={} \newtoks\RunTitle\RunTitle={}
\def\ShortTitle#1#2{\RunAuthor={#1}\RunTitle={#2}}
\headline={\ifnum\pageno=1{\hfil}\else
	\ifodd\pageno {\petit{\the\RunTitle}\hfil\folio}
	\else {\petit\folio\hfil{\the\RunAuthor}} \fi \fi}
\footline={\hfil}
\long\def\Title#1{\hrule\topglue3truecm
	\noindent{\bfbig#1}\vskip12pt\relax}
\long\def\Author#1#2{\noindent{\bf#1}\vskip6pt%
	\noindent{\petit#2}\vskip10pt%
	\noindent{\petit\today}\vskip32pt\relax}
\long\def\Thanks#1#2{{\parindent=20pt\baselineskip=9pt%
	\footnote{\nobreak${}^{#1}$}{\petit#2\par\vskip-9pt}}}
\long\def\Abstract#1{\noindent{\bf Summary.}\enspace#1\bigskip\relax}
\def\Acknowledgements{\goodbreak\vskip21pt\noindent{\bf Acknowledgements.}%
	\enspace\relax}


\def\PLabel#1{\xdef#1{\nobreak(p.\the\pageno)}}


\newcount\SECNO \SECNO=0
\newcount\SUBSECNO \SUBSECNO=0
\newcount\SUBSUBSECNO \SUBSUBSECNO=0
\def\Part#1{\SECNO=0\SUBSECNO=0\SUBSUBSECNO=0
	\vfill\eject\noindent{\partf Part #1}
	\nobreak\vskip12pt\noindent\kern0pt}
\def\Section#1{\SUBSECNO=0\SUBSUBSECNO=0 \advance\SECNO by 1
	\goodbreak\vskip21pt\leftline{\bf\the\SECNO .\ #1}
	\gdef\Label##1{\xdef##1{\nobreak\the\SECNO}}
	\nobreak\vskip12pt\noindent\kern0pt}
\def\SubSection#1{\SUBSUBSECNO=0 \advance\SUBSECNO by 1
	\goodbreak\vskip21pt\leftline{\it\the\SECNO.\the\SUBSECNO\ #1}
	\gdef\Label##1{\xdef##1{\nobreak\the\SECNO.\the\SUBSECNO}}
	\nobreak\vskip12pt\noindent\kern0pt}
\def\SubSubSection#1{\advance\SUBSUBSECNO by 1
	\goodbreak\vskip21pt\leftline{\rm\the\SECNO.\the\SUBSECNO.\the\SUBSUBSECNO\ #1}
	\gdef\Label##1{\xdef##1{\nobreak\the\SECNO.\the\SUBSECNO.\the\SUBSUBSECNO}}
	\nobreak\vskip12pt\noindent\kern0pt}


\long\def\Definition#1#2{\medbreak\noindent{\bf Definition%
	#1.\enspace}{\it#2}\medbreak\smallskip\relax}
\long\def\Theorem#1#2{\medbreak\noindent{\bf Theorem%
	#1.\enspace}{\it#2}\medbreak\smallskip\relax}
\long\def\Lemma#1#2{\medbreak\noindent{\bf Lemma%
	#1.\enspace}{\it#2}\medbreak\smallskip\relax}
\long\def\Proposition#1{\medbreak\noindent{\bf Proposition.%
	\enspace}{\it#1}\medbreak\smallskip\relax}
\long\def\Corollary#1#2{\medbreak\noindent{\bf Corollary%
	#1.\enspace}{\it#2}\medbreak\smallskip\relax}
\long\def\Conjecture#1#2{\medbreak\noindent{\bf Property%
	#1.\enspace}{\it#2}\medbreak\smallskip\relax}


\newcount\FOOTNO \FOOTNO=0
\long\def\Footnote#1{\global\advance\FOOTNO by 1
	{\parindent=20pt\baselineskip=9pt%
	\footnote{\nobreak${}^{\the\FOOTNO)}$}{\petit#1\par\vskip-9pt}%
	}\gdef\Label##1{\xdef##1{\nobreak\the\FOOTNO}}}


\newcount\EQNO \EQNO=0
\def\Eqno{\global\advance\EQNO by 1 \eqno(\the\EQNO)%
	\gdef\Label##1{\xdef##1{\nobreak(\the\EQNO)}}}


\newcount\FIGNO \FIGNO=0
\def\Fcaption#1{\global\advance\FIGNO by 1
	{\petit{\bf Fig. \the\FIGNO.~}#1}
	\gdef\Label##1{\xdef##1{\nobreak\the\FIGNO}}}
\def\Figure#1#2{\medskip\vbox{\centerline{\framebox{#1}}
	\centerline{\Fcaption{#2}}}\medskip\relax}

\newcount\TABNO \TABNO=0
\def\Tcaption#1{\global\advance\TABNO by 1
   {\petit{\bf Table. \the\TABNO.~}#1}
   \gdef\Label##1{\xdef##1{\nobreak\the\TABNO}}}
\def\Table#1#2{\medskip\vbox{\centerline{\Tcaption{#1}}\vskip4pt
   \centerline{#2}}\medskip\relax}

\newread\epsffilein    
\newif\ifepsffileok    
\newif\ifepsfbbfound   
\newif\ifepsfverbose   
\newif\ifepsfdraft     
\newdimen\epsfxsize    
\newdimen\epsfysize    
\newdimen\epsftsize    
\newdimen\epsfrsize    
\newdimen\epsftmp      
\newdimen\pspoints     
\pspoints=1bp          
\epsfxsize=0pt         
\epsfysize=0pt         
\def\epsfbox#1{\global\def\epsfllx{72}\global\def\epsflly{72}%
   \global\def\epsfurx{540}\global\def\epsfury{720}%
   \def\lbracket{[}\def\testit{#1}\ifx\testit\lbracket
   \let\next=\epsfgetlitbb\else\let\next=\epsfnormal\fi\next{#1}}%
\def\epsfgetlitbb#1#2 #3 #4 #5]#6{\epsfgrab #2 #3 #4 #5 .\\%
   \epsfsetgraph{#6}}%
\def\epsfnormal#1{\epsfgetbb{#1}\epsfsetgraph{#1}}%
\def\epsfgetbb#1{%
%
%
\openin\epsffilein=#1
\ifeof\epsffilein\errmessage{I couldn't open #1, will ignore it}\else
%
%
   {\epsffileoktrue \chardef\other=12
    \def\do##1{\catcode`##1=\other}\dospecials \catcode`\ =10
    \loop
       \read\epsffilein to \epsffileline
       \ifeof\epsffilein\epsffileokfalse\else
%
%
          \expandafter\epsfaux\epsffileline:. \\%
       \fi
   \ifepsffileok\repeat
   \ifepsfbbfound\else
    \ifepsfverbose\message{No bounding box comment in #1; using defaults}\fi\fi
   }\closein\epsffilein\fi}%
%
%
%
\def\epsfclipoff{\def\epsfclipstring{\ifepsfdraft\space clip\fi}}%
\epsfclipoff
\def\epsfsetgraph#1{%
   \epsfrsize=\epsfury\pspoints
   \advance\epsfrsize by-\epsflly\pspoints
   \epsftsize=\epsfurx\pspoints
   \advance\epsftsize by-\epsfllx\pspoints
%
%
   \epsfxsize\epsfsize\epsftsize\epsfrsize
   \ifnum\epsfxsize=0 \ifnum\epsfysize=0
      \epsfxsize=\epsftsize \epsfysize=\epsfrsize
      \epsfrsize=0pt
%
%
     \else\epsftmp=\epsftsize \divide\epsftmp\epsfrsize
       \epsfxsize=\epsfysize \multiply\epsfxsize\epsftmp
       \multiply\epsftmp\epsfrsize \advance\epsftsize-\epsftmp
       \epsftmp=\epsfysize
       \loop \advance\epsftsize\epsftsize \divide\epsftmp 2
       \ifnum\epsftmp>0
          \ifnum\epsftsize<\epsfrsize\else
             \advance\epsftsize-\epsfrsize \advance\epsfxsize\epsftmp \fi
       \repeat
       \epsfrsize=0pt
     \fi
   \else \ifnum\epsfysize=0
     \epsftmp=\epsfrsize \divide\epsftmp\epsftsize
     \epsfysize=\epsfxsize \multiply\epsfysize\epsftmp
     \multiply\epsftmp\epsftsize \advance\epsfrsize-\epsftmp
     \epsftmp=\epsfxsize
     \loop \advance\epsfrsize\epsfrsize \divide\epsftmp 2
     \ifnum\epsftmp>0
        \ifnum\epsfrsize<\epsftsize\else
           \advance\epsfrsize-\epsftsize \advance\epsfysize\epsftmp \fi
     \repeat
     \epsfrsize=0pt
    \else
     \epsfrsize=\epsfysize
    \fi
   \fi
%
%
   \ifepsfverbose\message{#1: width=\the\epsfxsize, height=\the\epsfysize}\fi
   \epsftmp=10\epsfxsize \divide\epsftmp\pspoints
   \vbox to\epsfysize{\vfil\hbox to\epsfxsize{%
      \ifnum\epsfrsize=0\relax
        \includegraphics{\ifepsfdraft}%
      \else
        \epsfrsize=10\epsfysize \divide\epsfrsize\pspoints
        \includegraphics{\ifepsfdraft}%
      \fi
      \hfil}}%
\global\epsfxsize=0pt\global\epsfysize=0pt}%
%
%
{\catcode`\%=12 \global\let\epsfpercent=
%
%
\long\def\epsfaux#1#2:#3\\{\ifx#1\epsfpercent
   \def\testit{#2}\ifx\testit\epsfbblit
      \epsfgrab #3 . . . \\%
      \epsffileokfalse
      \global\epsfbbfoundtrue
   \fi\else\ifx#1\par\else\epsffileokfalse\fi\fi}%
%
%
\def\epsfempty{}%
\def\epsfgrab #1 #2 #3 #4 #5\\{%
\global\def\epsfllx{#1}\ifx\epsfllx\epsfempty
      \epsfgrab #2 #3 #4 #5 .\\\else
   \global\def\epsflly{#2}%
   \global\def\epsfurx{#3}\global\def\epsfury{#4}\fi}%
%
%
\def\epsfsize#1#2{\epsfxsize}
%
%



\newcount\REFNO \REFNO=0
\newbox\REFBOX \setbox\REFBOX=\vbox{}
\def\BegRefs{\setbox\REFBOX\vbox\bgroup
	\parindent18pt\baselineskip9pt\petit}
\def\EndRefs{\par\egroup}
\def\Ref#1{\global\advance\REFNO by 1 \ifnum\REFNO>1\vskip3pt\fi
	\item{\the\REFNO .~}\xdef#1{\nobreak[\the\REFNO]}}
\def\References{\goodbreak\vskip21pt\leftline{\bf References}
	\nobreak\vskip12pt\unvbox\REFBOX\vskip21pt\relax}


\def\N{I\kern-.8ex N}
\def\Z{\raise.72ex\hbox{${}_{\not}$}\kern-.45ex {\rm Z}}
\def\Q{\raise.82ex\hbox{${}_/$}\kern-1.35ex Q} \def\R{I\kern-.8ex R}
\def\C{\raise.87ex\hbox{${}_/$}\kern-1.35ex C} \def\H{I\kern-.8ex H}

\def\D#1#2{{{\partial#1}\over{\partial#2}}}
\def\DV{{\rm D}\!{\rm V}}
\def\RP{{\R\!P}}



\Title{Orthogonal nets and Clifford algebras}
\Author{Alexander I.\ Bobenko and
	Udo J.\ Hertrich-Jeromin\Thanks{\ast}{Partially supported by
	the Alexander von Humboldt Stiftung and NSF grant DMS93-12087}}
	{Dept.\ Mathematics, Technical University Berlin,
	 D-10623 Berlin\newline
	 Dept.\ Math.\ \&\ Stat., GANG, University of Massachusetts,
	 Amherst, MA 01003}
\ShortTitle{A.Bobenko, U.Jeromin}
	{Orthogonal nets}
\Abstract{A Clifford algebra model for M\"obius geometry is presented.
	The notion of Ribaucour pairs of orthogonal systems in arbitrary
	dimensions is introduced, and the structure equations for adapted
	frames are derived. These equations are discretized and the geometry
	of the occuring discrete nets and sphere congruences is discussed
	in a conformal setting.
	This way, the notions of ``discrete Ribaucour congruences'' and
	``discrete Ribaucour pairs of orthogonal systems'' are obtained
	--- the latter as a generalization of discrete orthogonal systems
	in Euclidean space.
	The relation of a Cauchy problem for discrete orthogonal nets and
	a permutability theorem for the Ribaucour transformation of smooth
	orthogonal systems is discussed.}

\BegRefs
\Ref\Berger M.~Berger: {\it Geometry I\/}; Springer, Berlin 1987
\Ref\Blaschke W.~Blaschke: {\it Vorlesungen \"uber Differentialgeometrie
        III\/}; Springer, Berlin 1928
\Ref\DIso A.~Bobenko, U.~Pinkall: {\it Discrete isothermic surfaces\/};
        J.~reine angew.\ Math.\ {\bf 475} (1996) 187-208
\Ref\Bob A.~Bobenko: {\it Discrete Conformal Nets and Surfaces\/};
	GANG Preprint IV.27 (1996), to appear in P.~Clarkson, F.~Nijhoff
	(eds.), Proceedings SIDE II Conference, Canterbury, July 1996,
	Cambridge University Press, Cambridge
\Ref\Bobe A.~Bobenko: {\it Discrete Integrable Systems and Geometry\/};
	SFB 288 Preprint 298 (1997), to appear in the Proceedings of the
	International Congress of Mathematical Physics 1997, Brisbane,
	July 1997
\Ref\BoPi A.~Bobenko, U.~Pinkall: {\it Discretization of Surfaces and
        Integrable Systems\/}; SFB 288 Preprint 296 (1997),
        to appear in A.~Bobenko, R.~Seiler (eds.), {\it Discrete integrable
        Geometry and Physics\/}, Oxford Univ.\ Press, Oxford 1998
\Ref\CDS J.~Cie\'{s}li\'{n}ski, A.~Doliwa, P.~Santini: {\it The integrable
        Discrete Analogues of Orthogonal coordinate systems are
        multidimensional Circular lattices\/}; Phys.\ Lett.\ A {\bf 235}
	(1997) 480-488
\Ref\Darboux G.~Darboux: {\it Le\c{c}ons sur les syst\`{e}mes orthogonaux
        et les coordonn\'{e}es curvilignes\/}; Gauthier-Villars, Paris 1910
\Ref\DoSa A.~Doliwa, P.~Santini: {\it Geometry of Discrete Curves and
        Lattices and Integrable Difference equations\/}; Preprint (1997),
        to appear in A.~Bobenko, R.~Seiler (eds.), {\it Discrete integrable
        Geometry and Physics\/}, Oxford Univ.\ Press, Oxford 1998
\Ref\DMS A.~Doliwa, S.~Manakov, P.~Santini: {\it $\bar{\partial}$-reductions
	of the multidimensional quadrilateral lattice I: the multidimensional
	circular lattice\/}; Preprint Universit\`{a} di Roma (1997)
\Ref\Eisenhart L.~Eisenhart: {\it Transformations of Surfaces\/};
        Chelsea, New York 1962
\Ref\Fiechte W.~Fiechte: {\it Ebene M\"obiusgeometrie\/}; Manuscript 1979
\Ref\GaTsa E.~Ganzha, S.~Tsarev: {\it On superposition of the auto
        B\"acklund transformations for $(2+1)$-dimensional integrable
        systems\/}; Preprint (1996) solv-int 9606003
\Ref\Guichard C.~Guichard: {\it Sur les syst\`{e}mes triplement
        ind\'{e}termin\'{e}s et sur les syst\`{e}mes triple orthogonaux\/};
        Scientia 25, Gauthier-Villars, Paris 1905
\Ref\JTZ U.~Hertrich-Jeromin, E.~Tjaden, M.~Z\"urcher: {\it Cyclic
        systems and Guichard's nets\/}; MSRI preprint, Z\"urich 1996
\Ref\Suppl U.~Hertrich-Jeromin: {\it Supplement on Curved flats in the space
        of Point pairs and Isothermic surfaces: A Quaternionic calculus\/};
        Doc.\ Math.\ J.~DMV {\bf 2} (1997) 335-350
\Ref\JePe U.~Hertrich-Jeromin, F.~Pedit: {\it Remarks on the Darboux
        transform of isothermic surfaces\/};
        Doc.\ Math.\ J.~DMV {\bf 2} (1997) 313-333
\Ref\HJP U.~Hertrich-Jeromin, T.~Hoffmann, U.~Pinkall: {\it A discrete
        version of the Darboux transformation for isothermic surfaces\/};
	SFB 288 Preprint 239 (1996),
        to appear in A.~Bobenko, R.~Seiler (eds.), {\it Discrete integrable
        Geometry and Physics\/}, Oxford Univ.\ Press, Oxford 1998
\Ref\Infinitesimal U.~Hertrich-Jeromin: {\it The surfaces capable of
        division into infinitesimal Squares by their Curves of Curvature\/};
        GANG preprint, Amherst 1997
\Ref\KoSch B.~Konopelchenko, W.~Schief: {\it Three-dimensional integrable
        lattices in Euclidean spaces: Conjugacy and Orthogonality\/};
        Preprint (1997), to appear in Proc.\ Royal Soc.\ London
\Ref\Krichever I.~Krichever: {\it Algebraic-geometrical $n$-orthogonal
	curvilinear systems and solutions to the associativity equations\/};
	Preprint (1997)
\Ref\Salkowski E.~Salkowski: {\it Dreifach orthogonale Fl\"achensysteme\/};
        in Encyclopaedie der mathematischen Wissenschaften III.D 9,
        Teubner, Leipzig 1902
\Ref\Zakharov V.~Zakharov: {\it Description of the $n$-orthogonal curvilinear
        coordinate systems and Hamiltonian integrable systems of hydrodynamic
        type. Part I: Integration of the Lam\'{e} equations\/}; Preprint (1996),
	to appear in Duke Math.\ J.
\EndRefs


\Section{Introduction}
Triply orthogonal systems in Euclidean 3-space and Ribaucour sphere
congruences in the conformal 3-sphere were intensively studied around
the turn of the century --- among others by great geometers like Darboux
\Darboux, Bianchi, Guichard \Guichard\ (cf.\Salkowski), and Blaschke
\Blaschke. After this first wave of work, a calm period followed ---
until recently, when relations between ($n$-dimensional) orthogonal systems
and integrable system methods stimulated new interest. The corresponding
equations have been integrated by the $\bar{\partial}$-method \Zakharov\
and by the finite-gap integration theory \Krichever, and the Ribaucour
transformation for triply orthogonal systems have been studied in \GaTsa.
In differential geometry, there are relations with the theory of conformally
flat hypersurfaces (cf.\JTZ) and (the Darboux transformation) of isothermic
surfaces (cf.\Suppl,\Infinitesimal,\DIso,\JePe).
Circular nets as a natural geometrical discretization of triply orthogonal
coordinate systems have been introduced in \Bob\ and generalized for higher
dimensions in \CDS.
Geometric and analytic integrability of these nets were investigated in
\CDS, \DMS\ and \KoSch.
Also, discrete Ribaucour sphere congruences have been defined and studied
in \Bob, \HJP\ and \BoPi.

Both notions, the one of ``orthogonal systems'' as well as that of
``Ribaucour congruences'', are conformally invariant, i.e.\ invariant
under M\"obius transformations of the ambient space.
Nevertheless, in all the papers mentioned above (besides \JTZ), orthogonal
systems were treated in a Euclidean setting --- even though certain aspects
of the theory, e.g.\ the Ribaucour transformation for orthogonal systems,
seem to arise more naturally in a conformal setting.
Thus, we hope to make a small contribution in the growing flood of
publications in this field: the purpose of this paper is to introduce
a method that allows the treatment of smooth and discrete ``Ribaucour
pairs of orthogonal nets'' in M\"obius geometry\Footnote{In particular,
we also obtain a geometrical interpretation of the spinor Lax-representation
\DIso\ for isothermic surfaces.}.
Moreover, we expect this method not only to be advantageous in the theory
of discrete nets but also for the treatment of global questions in higher
dimensional M\"obius geometry (cf.\Suppl).

First, we give a short introduction to $n$-dimensional M\"obius geometry
--- in particular, we introduce a model using the Clifford algebra ${\cal C}$
of $(n+2)$-dimensional Minkowski space $\R^{n+2}_1$. To obtain geometric
interpretations for certain Clifford algebra elements, we use the fact that
the Clifford algebra ${\cal C}\cong\Lambda$ is isomorphic to the Grassmann
algebra as a vector space. On these M\"obius geometric objects, the spin
group ${\rm Spin}(\R^{n+2}_1)$ acts by M\"obius transformations.
Then, we discuss (smooth) orthogonal systems and Ribaucour congruences in
this setting. As a generalization, we introduce the notion of ``Ribaucour
pairs of orthogonal systems'' and derive the structure equations for
``adapted frames''.

These frame equations can easily be discretized: analyzing the occuring
discrete nets and ``discrete sphere congruences'', we obtain definitions
and characterizations for ``discrete Ribaucour congruences'' and ``discrete
Ribaucour pairs of orthogonal nets''. The latter generalizes the notion of
discrete orthogonal nets in Euclidean space (cf.\Bob,\CDS,\KoSch).
The occurence of certain (B\"acklund, Darboux, Ribaucour) transformations
is crucial in the relation of differential geometry and integrable system
theory (cf.\Bobe).
In the last part, we discuss the relation of a Cauchy problem for (Ribaucour
pairs of) discrete orthogonal systems (cf.\CDS,\DoSa) and a general version
of the permutability theorem for the Ribaucour transformation of (smooth)
orthogonal systems (cf.\Darboux,\GaTsa).
In fact, this permutability theorem shows how discrete orthogonal nets can
be obtained by repeatedly applying Ribaucour transformations to a given
orthogonal system.

\Section{M\"obius geometry and Clifford algebra}
In this section we will give a sketchy introduction to M\"obius
geometry --- especially, we intend to elaborate the Clifford algebra
model for M\"obius geometry and to relate it to the classical models.
Classically, the conformally compactified Euclidean space
$\R^n\cup\{p_{\infty}\}$ is considered the underlying space for the
group of M\"obius transformations \Blaschke\ --- those transformations
that map spheres (and planes --- that are considered as spheres
containing the point $p_{\infty}$ at infinity) to spheres.
It can be shown that the M\"obius group is generated by inversions
$$\matrix{p\mapsto m+{r^2\over|p-m|^2}(p-m)}$$ at spheres
(here, $m$ is the center and $r$ the radius of a sphere) and reflections
at planes. Via stereographic projection, $\R^n\cup\{p_{\infty}\}$ can be
identified with the $n$-sphere $S^n$ which is then embedded as an
absolute quadric in $\RP^{n+1}$. This allows to consider M\"obius
geometry as a subgeometry of projective geometry --- thus, to
linearize the M\"obius group: any M\"obius transformation of $S^n$
extends to a projective transformation of $\RP^{n+1}$ that fixes the
absolute quadric $S^n$.
For example, an inversion extends to a polar reflection of $\RP^{n+1}$.
In this model, hyperspheres $s\subset S^n$ are given as intersections
of projective hyperplanes with $S^n\subset\RP^{n+1}$. This way, they
can be identified with points in the ``outer space\Footnote{The
``inner space'' of $S^n\subset\RP^{n+1}$ can be defined as the set
of those points {\it not} lying on any (real) tangent line of $S^3$,
the outer space is its complement in $\RP^{n+1}\setminus S^n$.}''
$(\RP^{n+1})_o$ of $S^n$ via polarity --- a sphere $s\subset S^n$
is identified with the center of the cone that touches $S^n$ in $s$.

On the space of homogeneous coordinates of $\RP^{n+1}$ there is a
Lorentz scalar product --- unique up to scaling --- such that its
isotropic lines are exactly the points of $S^n$: light cone vectors
and unit vectors, $$\matrix{
  p\in L^{n+1}   &:=& \{v\in\R^{n+2}_1\,|\,\langle v,v \rangle=0\}\hfill\cr
  s\in S^{n+1}_1 &:=& \{v\in\R^{n+2}_1\,|\,\langle v,v \rangle=1\},\cr}$$
represent points and hyperspheres in $S^n$, respectively (cf.\Blaschke).
Incidence of a point and a hypersphere --- polarity in $\RP^{n+1}$ --- is
encoded as orthogonality: $p\in s$ if and only if $\langle p,s\rangle=0$.
And, a polar reflection in $\RP^{n+1}$, an inversion in $S^n$, is represented
by an ordinary orthogonal reflection $p\mapsto\pm(p-2\langle p,s\rangle s)$.
In this model, the metric subgeometries --- in particular Euclidean
geometry --- of M\"obius geometry are easily described: if $n_k\in\R^{n+2}_1$
is a vector with $k=-\langle n_k,n_k\rangle$, then the quadric
$Q_k:=\{p\in L^{n+1}\,|\,\langle p,n_k\rangle=1\}$ has constant sectional
curvature $k$; the stereographic projection $S^n\to\R^n$ becomes the central
projection $Q_1\to Q_0$ along the light cone generators.

Embedding the Minkowski space $\R^{n+2}_1$ into its Clifford algebra
${\cal C}$ additionally provides a useful algebraic structure.
As generators of ${\cal C}$, we consider a pseudo orthogonal basis
$(e_0,e_1,\dots,e_n,e_{\infty})$ of $\R^{n+2}_1$: $$\matrix{
  -e_i^2 &=& e_0e_{\infty}+e_{\infty}e_0 &=& 1\cr
  e_0^2 = e_{\infty}^2 &=& e_ie_j+e_je_i &=& 0\cr}$$
for $1\leq i\neq j\leq n$. Thus, the group ${\rm Spin}(\R^{n+2}_1)=
\{\Pi_{j=1}^{2m}s_j\,|\,s_j\in S^{n+1}_1\}$ acts as a double cover of the
isometry group $SO_1(n+2)$ on $\R^{n+2}_1$ via $p\mapsto\Phi^{-1}p\Phi$:
any orientation preserving M\"obius transformation is represented as
an (even) composition of inversions at spheres,
$$\matrix{p\mapsto\pm sps, & s\in S^{n+1}_1.\cr}$$
For our considerations, it will turn out crucial that the spin group and
its Lie algebra ${\fk spin}(\R^{n+2}_1)=\Lambda^2$ are both described as
subspaces of the same space ${\cal C}\cong\Lambda=\oplus_{i=0}^{n+2}\Lambda^i$
--- here, we identify the Clifford algebra ${\cal C}$ with the Grassmann
algebra over $\R^{n+2}_1$ as vector spaces. As the spin group acts
on $\Lambda^1$ via adjoint action, the Lie algebra acts on $\Lambda^1$
via the adjoint representation, $[v,\phi]=v\phi-\phi v$ for $v\in\Lambda^1$
and $\phi\in\Lambda^2$. In particular, for $0\leq i,j,k\leq\infty$, $i\neq j$,
$$[e_k,e_ie_j]=2(\langle e_j,e_k\rangle e_i-\langle e_i,e_k\rangle e_j).
  \Eqno\Label\EEE$$

The grading
${\cal C}\cong\Lambda=\oplus_{i=0}^{n+2}\Lambda^i$ of the Grassmann algebra
provides geometric interpretations for elements of the Clifford algebra:
clearly, points and hyperspheres in the conformal $n$-sphere $S^n$ are
represented by vectors $p,s\in\Lambda^1$ with $p^2=0$ and $s^2=-1$,
respectively. And, as an $m$-dimensional sphere $s$ in $S^n$ can be described
as the intersection of $n-m$ (orthogonal) hyperspheres $s_1,\dots,s_{n-m}\in
S^{n+1}_1\subset\Lambda^1$, we can identify it with the spacelike plane
${\rm span}\{s_1,\dots,s_{n-m}\}\subset\R^{n+2}_1$, i.e.\ $$\matrix{
  s=s_1\wedge\dots\wedge s_{n-m}\in\Lambda^{n-m} &{\rm with} &
 |s|^2:={\rm det}(\langle s_i,s_j\rangle)>0.\cr}$$
Thus, {\it pure\Footnote{Pure 2-vectors can be detected by the fact that
$s=s_1\wedge s_2$ for $s\in\Lambda^2$ if and only if $s\wedge s=0$, or
$s^2\in\R$. We use this fact in the table below.
Note that for pairwise orthogonal $s_1,\dots,s_k\in\Lambda^1$, the exterior
product coincides with the Clifford product, and $(\Pi_{i=1}^ks_i)^2=
(-1)^{k(k-1)\over2}\Pi_{i=1}^ks_i^2$.}} spacelike $(n-m)$-vectors
$s\in\Lambda^{n-m}$ represent $m$-spheres.
Choosing $\varepsilon:=(e_0-e_{\infty})\wedge e_1\wedge\dots\wedge e_n\wedge
(e_0+e_{\infty})\in\Lambda^{n+2}$, the linear isomorphism\Footnote{On non
isotropic pure $(n-m)$-vectors, this is Hodge duality, up to sign.}
$$\Lambda^{n-m}\ni v\mapsto\varepsilon v\in\Lambda^{m+2}$$
maps pure $(n-m)$-vectors to pure $(m+2)$-vectors --- reversing the signature
since $|\varepsilon v|^2=|\varepsilon|^2|v|^2=-|v|^2$. Consequently,
pure timelike $(m+2)$-vectors can be interpreted as $m$-spheres, too.
In particular,

\Lemma{}{$m+2$ points $p_i\in L^{n+1}$, $1\leq i\leq m+2$, in general
position determine a unique $m$-sphere $s$ containing all points, $p_i\in s$,
$$s=p_1\wedge\dots\wedge p_{m+2}\in\Lambda^{m+2}\cong\Lambda^{n-m},$$
and the $m+2$ points lie on a $(m-1)$-sphere if and only if
$p_1\wedge\dots\wedge p_{m+2}=0$.}

\noindent
In case four points $p_1,\dots,p_4\in\Lambda^1$ are not concircular, the
$\Lambda^4$-part of the quantity $p_1p_2p_3p_4+p_4p_3p_2p_1$ determines
a unique 2-sphere that contains the four points. Considering this 2-sphere
as a Riemann sphere, we can determine the (complex) cross ratio of the four
points:

\Lemma{}{Given four points $p_1,\dots,p_4\in\Lambda^1$, $p_i^2=0$, in the
conformal $n$-sphere their cross ratio is given by the quantity\Footnote{Note,
that $r$ does not depend on the choice of representative for the points $p_i$.}
$$\matrix{r=r_0+r_4
  ={p_1p_2p_3p_4+p_4p_3p_2p_1\over(p_1p_4+p_4p_1)(p_2p_3+p_3p_2)}
  \in\Lambda^0\oplus\Lambda^4}:\Eqno\Label\CrossRatio$$}

\noindent
up to complex conjugation, $r_0+i|r_4|$ is the complex cross ratio of the four
points on that sphere as a Riemann sphere --- since $r_0+i|r_4|$ is clearly
invariant under M\"obius transformations $p_i\mapsto\Phi^{-1}p_i\Phi$ it
suffices to check the formula using the identification $\C\ni z=x+iy\cong
|z|^2e_0+xe_1+ye_2+e_{\infty}\in L^{n+1}$, which is a straightforward
calculation. Moreover,

\Corollary{}{Four points $p_1,\dots,p_4\in L^{n+1}$ in the conformal $n$-sphere
are concircular if and only if their cross ratio \CrossRatio\ is real i.e.\
$r_4=0$.}

\noindent
For reference, we have summarized the essentials of the previous discussions
in Table 1 where, for simplicity of notation, we restrict to the case $n=3$.

\Section{Ribaucour pairs of orthogonal systems}
In this section, we are going to discuss briefly the geometry of smooth
orthogonal systems and Ribaucour congruences (cf.\Salkowski,\JTZ) ---
in particular, we will derive suitable frame equations.
Later, we will use discretizations of these frame equations to obtain
definitions for discrete analogs of smooth orthogonal systems and
Ribaucour congruences --- and, more general, of Ribaucour pairs of
orthogonal systems.

\Definition{~(orthogonal system)}{A system of $n$ 1-parameter families
of hypersurfaces in $\R^n$ is called an orthogonal system if any two
hypersurfaces from different families intersect orthogonally.}

\noindent
Obviously, the notion of an orthogonal system is invariant under M\"obius
transformations of the ambient space $\R^n$; or, more general, it is invariant
under conformal changes of the ambient space's metric.
Consequently, we will focus on conformally invariant properties of orthogonal
systems --- even though, we will do some calculations in a Euclidean setting
where it seems more convenient.

Parametrizing the $n$ families of orthogonal hypersurfaces one obtains an
orthogonal coordinate system $(t_1,\dots,t_n):V\subset\R^n\to\R^n$ resp.\
a parametrization $f:U\subset\R^n\to\R^n\subset L^{n+1}$ with $\D{}{t_i}f
\perp\D{}{t_j}f$ for any pair $1\leq i\neq j\leq n$.
Thus, given an orthogonal system $f$, we can introduce an {\it adapted framing}
$\Phi:U\to{\rm Spin}(\R^{n+2}_1)$ with
$$\matrix{\hfill f&=&\kern1.35em\Phi^{-1}e_0\Phi\hfill&\cr
  \hfill \D{}{t_i}f=:f_i&=&l_i\cdot\Phi^{-1}e_i\Phi,\hfill&1\leq i\leq n\cr}
  \Eqno\Label\OSFrame$$
where $(e_0,e_1,\dots,e_n,e_{\infty})$ denotes any pseudo orthonormal basis
of $\R^{n+2}_1$ with $\R^n\cong\{y\in L^{n+1}\,|\,ye_{\infty}+e_{\infty}y=1\}$
and $l_i:U\to\R$ are {\it Lam\'{e}'s functions} \Salkowski.
In terms of the adapted frame $\Phi$, the fact that $f:U\to\R^n$ takes values
in Euclidean space is encoded by $e_{\infty}=\Phi^{-1}e_{\infty}\Phi$.

In our investigations, a key role will be played by the following

\Theorem{~(Dupin)}{In an orthogonal system in $\R^n$, the intersection of
$n-1$ hypersurfaces from different families is a curvature line for any of
the intersecting hypersurfaces.}


\vfill\eject\topglue2truecm
\Table{M\"obius geometry in different models}{%
\long\def\entry#1#2{\matrix{\hbox{\petit#1}\hfill\cr#2\hfill\cr}\hfill}
$\matrix{\cr
\R^3\cong\hbox{Im}\H \hfill &
S^3\subset\RP^4 \hfill &
L^4\subset\R^5_1 \hfill &
{\cal C}\cong\Lambda=\oplus_{i=0}^5\Lambda^i \hfill \cr
\noalign{\vskip5pt\hrule\vskip5pt}\cr
\entry{point:}{p\in\R^3} &
\entry{point on $S^3$:}{p\in S^3} &
\entry{lightlike vector:}{p\in L^4} &
\entry{nilpotent vector:}{p\in\Lambda^1, p^2=0: L^4} \cr
\noalign{\vskip7pt}\cr
\entry{sphere:}{s\subset\R^3} &
\entry{point outside $S^3$:}{s\in(\RP^4)_o} &
\entry{unit vector:}{s\in S^4_1} &
\entry{anti-involutive vector:}{s\in\Lambda^1, s^2=-1: S^4_1} \cr
\noalign{\vskip0pt}\cr
\entry{plane:}{t\subset\R^3} &
\entry{\dots containing $p_{\infty}$:}{t\in T_{p_{\infty}}S^3} &
\entry{\dots perpendicular to $p_{\infty}$:}{t\in S^4_1,
	\langle t,p_{\infty}\rangle=0} &
\entry{\dots anti-symmetric to $p_{\infty}$:}{t\in S^4_1,
	{tp_{\infty}+p_{\infty}t\over2}=0} \cr
\noalign{\vskip7pt}\cr
\entry{circle:}{c\subset\R^3} &
\entry{line outside $S^3$:}{c\subset(\RP^4)_o} &
\entry{spacelike plane:}{c\in G_+(2,3)=:C^6_2} &
\entry{anti-involutive bivector:}{c\in\Lambda^2, c^2=-1: C^6_2} \cr
\noalign{\vskip0pt}\cr
\entry{line:}{l\subset\R^3} &
\entry{\dots containing $p_{\infty}$:}{l\subset T_{p_{\infty}}S^3} &
\entry{\dots perpendicular to $p_{\infty}$:}{l\in C^6_2,
	\langle l,p_{\infty}\rangle=0} &
\entry{\dots symmetric to $p_{\infty}$:}{l\in C^6_2,
	{lp_{\infty}-p_{\infty}l\over2}=0} \cr
\noalign{\vskip5pt\hrule\vskip5pt}\cr
\entry{incidence:}{p\in s} &
\entry{polarity:}{s\in T_pS^3} &
\entry{orthogonality:}{\langle p,s\rangle=0} &
\entry{anti-symmetry:}{ps+sp=0} \cr
\noalign{\vskip7pt}\cr
\entry{incidence:}{p\in c} &
\entry{polarity:}{c\subset T_pS^3} &
\entry{orthogonality:}{\langle p,s\rangle=0} &
\entry{symmetry:}{pc-cp=0} \cr
\noalign{\vskip7pt}\cr
\entry{intersection:}{s_1\cap_{\perp}s_2} &
\entry{polarity:}{s_1\in\hbox{Pol}[s_2]} &
\entry{orthogonality:}{\langle s_1,s_2\rangle=0} &
\entry{anti-symmetry:}{s_1s_2+s_2s_1=0} \cr
\noalign{\vskip0pt}\cr
\entry{\dots angle:}{s_1\cap_{\varphi}s_2} &
&
\entry{scalar product:}{\langle s_1,s_2\rangle=\cos\varphi} &
\entry{real part:}{{s_1s_2+s_2s_1\over2}=\cos\varphi} \cr
\noalign{\vskip7pt}\cr
\entry{intersection:}{c\cap_{\perp}s} &
\entry{polarity:}{c\subset\hbox{Pol}[s]} &
\entry{orthogonality:}{\langle s,c\rangle=0} &
\entry{symmetry:}{cs-sc=0} \cr
\noalign{\vskip7pt}\cr
\entry{cross ratio:}{{(p_1-p_2)\over(p_2-p_3)}{(p_3-p_4)\over(p_4-p_1)}} &
&
\entry{cross ratio ($l_{ij}=\langle p_i,p_j\rangle$):}
	{{{l_{12}l_{34}+l_{14}l_{23}-l_{13}l_{24}+\sqrt{|(l_{ij})|}}
	\over{2l_{14}l_{23}}}} &
\entry{cross ratio:}{{{p_1p_2p_3p_4+p_4p_3p_2p_1}
	\over{(p_1p_4+p_4p_1)(p_2p_3+p_3p_2)}}} \cr
\noalign{\vskip5pt\hrule\vskip5pt}\cr
\entry{inversion:}{p\mapsto m-{r^2\over p-m}} &
\entry{polar reflection}{\dots} &
\entry{reflection:}{p\mapsto\pm(p-2\langle p,s\rangle s)} &
\entry{reflection:}{p\mapsto\pm sps} \cr
}$}\Label\Models
\vfill\eject

\noindent
Consequently, restricting a parametrization $f:U\to\R^n$ of an orthogonal
system to any of the (coordinate) hypersurfaces $t_i=const$ provides a
curvature line parametrization of the hypersurface\Footnote{Note, that
higher dimensional hypersurfaces usually do not carry curvature line
coordinates.}.
Moreover, an adapted framing \OSFrame\ restricts to a principal framing
of any of the hypersurfaces of an orthogonal system:
$n_i:=\Phi^{-1}e_i\Phi$ yields a unit normal field for the hypersurface
$t_i=const$ and the remaining directions $\Phi^{-1}e_j\Phi$ are its principal
directions with curvatures $\kappa_{ij}=-{1\over l_il_j}\D{}{t_i}l_j$, i.e\
$[e_i,\Phi_j\Phi^{-1}]=-\kappa_{ij}l_je_j$ where $\Phi_j:=\D{}{t_j}\Phi$.
Thus (cf.\EEE), $$\matrix{\Phi_j&=&
  l_j(e_{\infty}+{1\over2}\sum_{i=1}^n\kappa_{ij}e_i)e_j\cdot\Phi}.
  \Eqno\Label\EuclFrame$$
In terms of a parametrization $f:U\to L^{n+1}$ of an orthogonal system,
a conformal change of the ambient space's metric can be modeled by a
``conformal deformation'' $f\mapsto e^uf$.
For an adapted frame \OSFrame, this has the effect that the second point of
intersection of the spheres $\Phi^{-1}e_j\Phi$  is generally
not constant any more, i.e.\ that these spheres cannot be interpreted as
planes in a Euclidean space any longer: $[e_{\infty},\Phi_j\Phi^{-1}]=
{1\over2}l_j\sum_{i=1}^n\sigma_{ij}e_i$ with the (symmetric) Schouten
tensor $(\sigma_{ij})$ of the induced metric $\sum_{i=1}^nl_i^2dt_i$.
Now, $$\matrix{
  \Phi_j &=& l_j[(e_{\infty}+{1\over2}\sum_{i=1}^n\gamma_{ij}e_i)e_j
  +e_0 ({1\over2}\sum_{i=1}^n\sigma_{ij}e_i)]\cdot\Phi}.\Eqno\Label\ConfFrame$$

The other key concept we are going to examine is the one of a

\Definition{~(Ribaucour congruence)}{A hypersurface is said to envelope a
sphere congruence (an $(n-1)$-parameter family of spheres) if, at every point,
the hypersurface has first order contact with a sphere of the congruence,
and\newline
a sphere congruence is said to be Ribaucour if the curvature lines on its
two envelopes do correspond.}

\noindent
A more technical description for the envelopes of a sphere congruence and
for Ribaucour congruences will prove useful (cf.\Blaschke,\JTZ):

\Lemma{}{$f:M^{n-1}\to L^{n+1}$ envelopes a sphere congruence $s:M^{n-1}\to
S^{n+1}_1$ if and only if $\langle s,f\rangle=0$ and $\langle ds,f\rangle=0$}

\noindent
Thus, the two envelopes $f,\hat{f}:M^{n-1}\to L^{n+1}$ of a (regular) sphere
congruence $s:M^{n-1}\to S^{n+1}_1$ can be interpreted as its two isotropic
normal fields --- that are uniquely determined up to rescalings.
And, since the principal curvature directions of the two envelopes coincide
with the principal directions of $s$ with respect to $f$ and $\hat{f}$,
respectively, as normal fields,

\Lemma{}{$s:M^{n-1}\to S^{n+1}_1$ is a Ribaucour congruence if and only if
its normal bundle is flat, i.e.\ if $d\langle df,\hat{f}\rangle=0$ for its
two envelopes $f,\hat{f}:M^{n-1}\to L^{n+1}$.}

\noindent
Consequently, the two isotropic normal fields of a Ribaucour congruence
--- its envelopes --- can be normalized to be {\it parallel} sections of
the normal bundle.
If we additionally {\it assume} the existence of principal curvature line
coordinates we obtain two orthogonal systems of codimension 1 in the
conformal $n$-sphere, with $\D{}{t_i}f/\!/\D{}{t_i}\hat{f}/\!/\D{}{t_i}s$.
Thus, an {\it adapted frame} $\Phi:U\subset\R^{n-1}\to{\rm Spin}(\R^{n+2}_1)$,
$$\matrix{
  \Phi^{-1}e_0\Phi=f,&\Phi^{-1}e_n\Phi=s,&\Phi^{-1}e_{\infty}\Phi=\hat{f},\cr
  {\rm and} & \Phi^{-1}e_i\Phi/\!/\D{}{t_i}f,&1\leq i\leq n-1, \cr}
  \Eqno\Label\RFrame$$
can be introduced.
Again, Dupin's theorem holds, i.e.\ $[e_i,\Phi_j\Phi^{-1}]\wedge e_j=0$
for $0\leq i\leq\infty$ and $1\leq j\leq n-1$, $j\neq i$: and consequently,
we find $$\matrix{ \Phi_j &=& 
  (l_je_{\infty}+{1\over2}\sum_{i=1}^{n-1}\gamma_{ij}e_i+{1\over2}a_je_n
  +\hat{l}_je_0)e_j\cdot\Phi}\Eqno\Label\RibFrame$$
with suitable functions $l_j,\hat{l}_j,a_j,\gamma_{ij}:U\to\R$,
$\gamma_{jj}=0$.
A rather remarkable property of such a ``Ribaucour pair of orthogonal
systems'' motivates a definition in case of arbitrary codimension:
any two corresponding subnets $t_i=const$ of $f$ and $\hat{f}$ do not only
envelope the sphere congruence $s$ but also $s_i:=\Phi^{-1}e_i\Phi$ since
$s_i\perp\D{}{t_j}f,\D{}{t_j}\hat{f}$ for $j\neq i$ --- and hence, both
subnets envelope the congruence $s\cdot s_i:U\to\Lambda^2$ of $(n-2)$-spheres.
Moreover, $\D{}{t_j}s\perp s_i$ for $j\neq i$ which shows that both subnets
carry well defined curvature lines (independent of the normal direction)
--- that, as before, do correspond on both envelopes and are given by the
coordinate directions $t_j$.
Clearly, these facts hold true for any pair of corresponding subnets of any
dimension.
Thus,

\Definition{}{A map $f:M^m\to S^n$ into the conformal $n$-sphere is said
to envelope a congruence of $m$-spheres $s:M^m\to\Lambda^{n-m}$ if, at
every point $p\in M^m$, $f$ has first order contact with the corresponding
sphere $s(p)$;\newline
a congruence of $m$-spheres $s:M^m\to\Lambda^{n-m}$ is called Ribaucour if
it has two envelopes with well defined curvature lines that do correspond;
and\newline
two nets $f,\hat{f}:U\subset\R^m\to S^n$ in the conformal $n$-sphere are
said to form a Ribaucour pair of orthogonal systems if any corresponding
$k$-dimensional subnets, $1\leq k\leq m$, envelope a Ribaucour congruence
of $k$-spheres.}

\noindent
In particular, the coordinate lines in the nets of a Ribaucour pair of
orthogonal systems intersect pairwise orthogonal since they arise as
curvature lines --- in case $m=n$ as curvature lines of lower dimensional
subnets.
We also explicitely allow one of the nets, $f$ or $\hat{f}$, to degenerate:
this way, the classical orthogonal systems in Euclidean $n$-space (coupled
with the point at infinity) appear as special cases of Ribaucour pairs of
orthogonal systems.

Finally, we want to characterize adapted frames for Ribaucour pairs of
orthogonal systems: we call $\Phi:U\subset\R^m\to{\rm Spin}(\R^{n+2}_1)$
an {\it adapted frame} if $$\matrix{
  f = \Phi^{-1}e_0\Phi, & \hat{f} = \Phi^{-1}e_{\infty}\Phi,\cr}$$
and if $s_i:=\Phi^{-1}e_i\Phi$ are the principal directions of $f$ and
$\hat{f}$ for $1\leq i\leq m$ and form a parallel orthonormal basis
of $\{f,\hat{f},df(TU)\}^{\perp}$ for $m+1\leq i\leq n$.
Then,

\Proposition{$\Phi:U\subset\R^m\to{\rm Spin}(\R^{n+2}_1)$ is an
adapted frame for a Ribaucour pair $f,\hat{f}:U\to L^{n+1}$ of
orthogonal systems in the conformal $n$-sphere iff
$$\matrix{ \Phi_j &=& 
  (l_je_{\infty}
  +{1\over2}\sum_{i=1}^m\gamma_{ij}e_i
  +{1\over2}\sum_{i=1}^{n-m}a_{ij}e_i
  +\hat{l}_je_0
  )e_j\cdot\Phi}\Eqno\Label\ROSFrame$$
for $1\leq j\leq m$, and suitable functions
$l_i,\hat{l}_i,\gamma_{ij},a_{ij}:U\to\R$, $\gamma_{ii}=0$.}

\noindent
At this point, we are prepared to formulate the discrete frame
equations that we will use to define (algebraically) the discrete
analogs for Ribaucour pairs of orthogonal systems:

\Section{Discrete orthogonal systems and Ribaucour congruences}
The basic idea is to discretize the frame equations we derived in the
previous section, i.e.\ to discretize the first order Taylor expansion
of the adapted frames $\Phi:U\to{\rm Spin}(\R^{n+2}_1)$:
for infinitesimal $0\simeq\varepsilon\in\R$, $t\in U$ and the standard
direction vectors $t_j=(\delta_{1j},\dots,\delta_{nj})$, we have
$$\Phi(t+\varepsilon t_j)\simeq\Phi(t)+\Phi_j(t)\cdot\varepsilon
  =[1+\varepsilon\Phi_j\Phi^{-1}(t)]\cdot\Phi(t).$$
A discrete analog of this equation can easily be formulated if $1$ and
$\Phi_j\Phi^{-1}$ can be combined to take values in ${\rm Spin}(\R^{n+2}_1)$
--- since $\Phi(t+\varepsilon t_j)$ and $\Phi(t)$ lie in ${\rm Spin}
(\R^{n+2}_1)$ the product $\Phi(t+\varepsilon t_j)\Phi^{-1}(t)\simeq[1+
\varepsilon\Phi_j\Phi^{-1}(t)]$ should, too.
Thus, if $[\varepsilon_{j1}+\varepsilon_{j2}\Phi_j\Phi^{-1}(t)]\in
{\rm Spin}(\R^{n+2}_1)$ for some {\it finite} numbers $\varepsilon_{j1},
\varepsilon_{j2}\in\R$ we can just use the {\it sructure} of
$[\varepsilon_{j1}+\varepsilon_{j2}\Phi_j\Phi^{-1}(t)]$ to define the
discrete analogs of the frame equations.

Having this in mind, an examination of the frame equations \EuclFrame\
for an adapted frame of an orthogonal system in {\it Euclidean ambient
space} $\R^n$ as well as \RibFrame\ for an adapted frame of a Ribaucour
congruence in the conformal $n$-sphere --- or, more general, the frame
equations \ROSFrame\ for a Ribaucour pair of orthogonal systems in the
conformal $n$-sphere --- can be discretized following that approach:
here, {\it any} combination $[1+\varepsilon\Phi_j\Phi^{-1}]$ with
$\varepsilon\in\R$ is a {\it pure} bivector --- and, consequently,
a suitable normalization takes values in ${\rm Spin}(\R^{n+2}_1)$.
The situation is different in the case of the adapted frame equations
\ConfFrame\ for a single orthogonal system in the conformal $n$-sphere:
in that case, it will be hard to control whether a combination
$1+\varepsilon\Phi_j\Phi^{-1}$ is a pure bivector, and therefore, whether
a suitable normalization can ever take values in ${\rm Spin}(\R^{n+2}_1)$.

Thus, we will stick to the case of Ribaucour pairs of orthogonal systems:
here, the structure equations for a discrete frame
$\Phi:\Gamma\subset\Z^m\to{\rm Spin}(\R^{n+2}_1)$ will be of the form
$$\matrix{\Phi(t+t_j)=e_js_j(t+{1\over2}t_j)\cdot\Phi(t)}
  \Eqno\Label\DiscreteFrame$$
with suitable vector functions $s_j:\Gamma_j\to S^{n+1}_1$ --- where we use
the notation $t+{1\over2}t_j$ for the edges in $j$-direction of $\Gamma$,
and $\Gamma_j$ for the lattice formed by these edges.
With this ansatz for $\Phi$, we examine the geometry of the two maps
$$\matrix{F:=\Phi^{-1}e_0\Phi:\ \Gamma\to L^{n+1}\cr
    \hat{F}:=\Phi^{-1}e_{\infty}\Phi:\Gamma\to L^{n+1}\cr}
    \Eqno\Label\Maps$$
--- which will lead us to the definition of {\it discrete Ribaucour pairs
of orthogonal systems}:
first, we notice that two points $F(t)$ and $F(t+t_j)$ lie symmetric with
respect to the sphere\Footnote{Note, that the second equality shows that
$S_j(t+{1\over2}t_j)$ depends symmetrically on the endpoints $t$ and $t+t_j$
of the edge $t+{1\over2}t_j$.}
$$\matrix{S_j(t+{1\over2}t_j):= \Phi^{-1}(t)e_j\Phi(t+t_j)
  = \Phi^{-1}(t+t_j)e_j\Phi(t)}:\Eqno\Label\ElementSphere$$
$$\matrix{F(t+t_j)=S_j(t+{1\over2}t_j)F(t)S_j(t+{1\over2}t_j)
  },$$
i.e.\ $F(t+t_j)$ is obtained from $F(t)$ by inversion at the sphere
$S_j(t+{1\over2}t_j)$.
Since the same holds true for the net $\hat{F}$ we deduce that any four
points $F(t)$, $\hat{F}(t)$, $F(t+t_j)$ and $\hat{F}(t+t_j)$ lie on a circle
--- as a calculation of their cross ratio \CrossRatio $$\matrix{
  \DV[\hat{F}(t),F(t),F(t+t_j),\hat{F}(t+t_j)]=-{1\over{
  \langle e_0,s_j(t+{1\over2}t_j)\rangle
  \langle e_{\infty},s_j(t+{1\over2}t_j)\rangle}}}$$ confirms
--- and the point pairs $(F(t),F(t+t_j))$ and $(\hat{F}(t),\hat{F}(t+t_j))$
do not seperate each other on the circle.
Conversely, if two corresponding edges of $F$ and $\hat{F}$ do not both
degenerate at the same time, $F(t)\neq F(t+t_j)$ or $\hat{F}(t)\neq
\hat{F}(t+t_j)$, then the corresponding sphere $S_j(t+{1\over2}t_j)$ can
be geometrically (i.e.\ up to sign) reconstructed from the four points
--- since they are concircular and do not seperate.
As another simple consequence we find that if one of the maps $F$ or $\hat{F}$
is constant, say $\hat{F}\equiv e_{\infty}$, then all the spheres $S_j$,
$1\leq j\leq m$, are planes, i.e.\ they contain $e_{\infty}$, the ``point
at infinity'': $S_je_{\infty}+e_{\infty}S_j=0$. In that case,
$F:\Gamma\to\R^n=\{y\in L^{n+1}\,|\,e_{\infty}y+ye_{\infty}=1\}$ takes values
in Euclidean space.

As in smooth differential geometry, a key tool in understanding the geometry
of the nets $F(\Gamma)$ and $\hat{F}(\Gamma)$ are the integrability conditions
$$\matrix{e_js_j(t+t_i+{1\over2}t_j)\cdot e_is_i(t+{1\over2}t_i)
  = e_is_i(t+t_j+{1\over2}t_i)\cdot e_js_j(t+{1\over2}t_j)}\Eqno\Label\MCs$$
of the structure equations \DiscreteFrame.
Using these ``Maurer-Cartan equations'', it is a straightforward calculation
to obtain the cross ratio \CrossRatio\ of an elementary quadrilateral with
vertices $F(t)$, $F(t+t_i)$, $F(t+t_i+t_j)$ and $F(t+t_j)$
$$\matrix{\DV[F(t),F(t+t_i),F(t+t_i+t_j),F(t+t_j)]=-
  {{\langle e_0,s_i(t+{1\over2}t_i)\rangle
    \langle e_0,s_i(t+t_j+{1\over2}t_i)\rangle}\over
   {\langle e_0,s_j(t+{1\over2}t_j)\rangle
    \langle e_0,s_j(t+t_i+{1\over2}t_j)\rangle}
  }}.$$
Since the cross ratio of four points is real if and only if the four points
lie on a circle we find that not only corresponding point pairs of $F$ and
$\hat{F}$ are concircular:

\Lemma{}{The vertices of any elementary quadrilateral of the net $F$, or
$\hat{F}$, as well as the endpoints of corresponding edges of $F$ and
$\hat{F}$ are concircular.}

\noindent
Together with the investigations in \BoPi\ or \Infinitesimal\ and Dupin's
theorem, this might motivate the following

\Corollary{~and Definition}{The two nets $F,\hat{F}:\Gamma^m\to S^n$ in
the conformal $n$-sphere are  discrete curvature line nets, i.e.\ the
vertices of any elementary quadrilateral (2-cell) have real cross ratio.}

\noindent
The above lemma directly generalizes to higher dimension: for example, we may
consider corresponding elementary quadrilaterals (2-cells) of the nets $F$
and $\hat{F}$.
Then, by the lemma, the vertices of the elementary quadrilateral of $F$ lie
on a circle, as well as the two endpoints of an edge and their corresponding
points of the net $\hat{F}$ do. These two circles are contained in a 2-sphere.
Now, taking another pair of edges that shares one pair of endpoints with the
former pair of edges, another circle is obtained. Since three of the points on
this circle already lie on the 2-sphere, the whole circle does. Consequently,
by symmetry, all the vertices of the elementary 2-cell of $\hat{F}$ also lie
on that 2-sphere.
A similar argument shows that the vertices of an elementary 3-cell of one of
the nets, $F$ or $\hat{F}$, lie on a 2-sphere --- and, finally, the same
arguments also prove similar statements for higher dimensional elementary
cells.
Thus, as a generalization of the previous lemma, we obtain the following

\Lemma{}{The vertices of any elementary $k$-cell of the net $F$, or $\hat{F}$,
lie on a $(k-1)$-sphere, and the vertices of any corresponding $k$-cells of
$F$ and $\hat{F}$ lie on a $k$-sphere, for $1\leq k\leq m$.}

\noindent
This lemma motivates the definition of a discrete Ribaucour congruence 
and its envelopes (cf.~\JTZ, \HJP), and of a discrete Ribaucour pair of
orthogonal systems:

\Corollary{~and Definition}{The two nets $F,\hat{F}:\Gamma^m\to S^n$ in
the conformal $n$-sphere envelope a discrete Ribaucour sphere congruence%
\Footnote{Note, that the spheres $S$ only depend on the elementary $m$-cells
of $\Gamma^m$, i.e.\ the domain of $S$ is the dual lattice $\Gamma^{\ast}$
of $\Gamma$.}
$$S=\Phi^{-1}e_0\wedge s_1\wedge\dots\wedge s_m\wedge e_{\infty}\Phi:
  \Gamma^{\ast}\to\Lambda^{m+2}\cong\Lambda^{n-m},\Eqno\Label\Congruence$$
i.e.\ all corresponding elementary $m$-cells of the two nets $F$ and $\hat{F}$
lie on $m$-spheres $S(t+{1\over2}(1,\dots,1))$.\newline
Moreover, the nets form a Ribaucour pair of orthogonal nets, i.e.\ all
corresponding elementary $k$-cells, $1\leq k\leq m$, lie on $k$-spheres
and the points of corresponding 1-cells (edges) do not seperate on the
circles they lie on.}

\noindent
A discrete map $S:\Gamma^{\ast}\to\Lambda^{n-m}$ into the space of $m$-spheres
in the conformal $n$-sphere generally has no envelopes: for the notion of
``envelopes'' to make sense one has to require the spheres
$S(t+{1\over2}(\varepsilon_1,\dots,\varepsilon_m))$, $t\in\Gamma$ and
$\varepsilon_j\in\{\pm1\}$, of each elementary $m$-cell in $\Gamma^{\ast}$
to intersect in a point pair.
In that case, it makes sense to speak of a ``discrete sphere congruence''.
Note, that for the notion of a ``discrete Ribaucour sphere congruence''
we {\it assume} the existence of two enveloping discrete curvature line nets.
Thus, the two envelopes can be reconstructed from a discrete Ribaucour
congruence: since we additionally know that the endpoints of corresponding
edges of the two envelopes should not seperate on the circles they lie on
(cf.\ElementSphere) there is a unique ordering on the pairs of corresponding
points.

From the previous lemma, it also becomes clear that --- analogous to the
smooth case --- any corresponding subnets of a discrete Ribaucour pair of
orthogonal nets form themselves a discrete Ribaucour pair.

With the equations \ElementSphere, the frame $\Phi$ can as well be
``integrated'' from the spheres $S_j$, $1\leq j\leq m$ --- and, since
the spheres $S_j$ can be (up to sign) constructed from the nets $F$
and $\hat{F}$,

\Theorem{}{All three, the discrete frame $\Phi:\Gamma\to{\rm Spin}(\R^{n+2}_1)$
satisfying \DiscreteFrame, the corresponding discrete Ribaucour pair of
orthogonal nets $F,\hat{F}:\Gamma\to L^{n+1}$, and the system of element
spheres $S_j:\Gamma_j\to S^{n+1}_1$ ($1\leq j\leq m$) can be constructed
from one of them.
If $m<n$, they can also be reconstructed from the discrete Ribaucour
congruence $S:\Gamma^{\ast}\to\Lambda^{n-m}$ of $m$-spheres enveloped
by $F$ and $\hat{F}$.}

\noindent
In these constructions, the frame is only determined up to a pointwise
sign change. In terms of the element spheres $S_j$, such sign changes
are reflected by a sign ambiguity, too. However, the spheres $S_j$ have
to satisfy a integrability condition that is equivalent to the Maurer-Cartan
equations \MCs\ and which restricts this sign ambiguity for the $S_j$:

\Lemma{~(Maurer-Cartan equations)}{Given maps $S_j:\Gamma_j\to S^{n+1}_1$
there exists a frame $\Phi:\Gamma\to{\rm Spin}(\R^{n+2}_1)$ satisfying
\ElementSphere\ for $1\leq j\leq m$ if and only if
$$\matrix{0=S_i(t+{1\over2}t_i)S_j(t+t_i+{1\over2}t_j)
  +S_j(t+{1\over2}t_j)S_i(t+t_j+{1\over2}t_i)}.\Eqno\Label\MCS$$}

\noindent
It follows that any four spheres $S_i(t+{1\over2}t_i)$, $S_j(t+{1\over2}t_j)$,
$S_j(t+t_i+{1\over2}t_j)$ and $S_i(t+t_j+{1\over2}t_i)$ intersect in a common
$(n-2)$-sphere. Thus, interpreted as points in projective $(n+1)$-space (cf.\
Table \Models), the four spheres are collinear.
Again, this picture can be directly generalized to higher dimensions:
\Figure{\hbox{\epsfbox{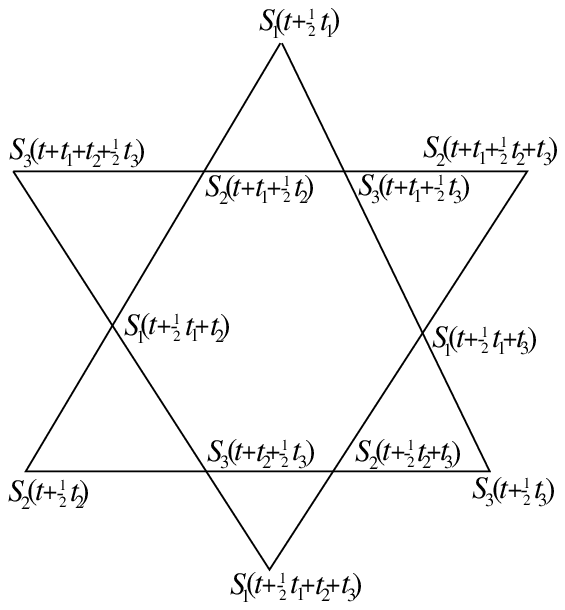}\vrule\epsfbox{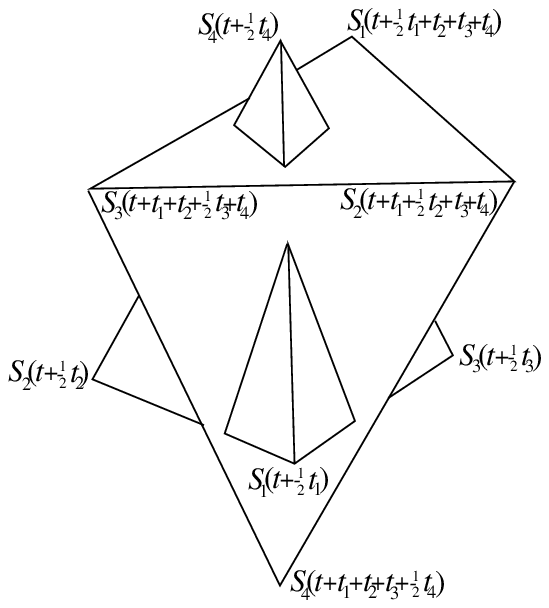}}}{Images
  of elementary 3- and 4-cells of $F$ in $\RP^{n+1}$}\Label\KCells\noindent
for $k\geq2$, any two elementary $k$-cells that belong to an elementary
$(k+1)$-cell share an elementary $(k-1)$-cell. Thus, assuming that all
spheres belonging to an elementary $k$-cell are contained in a $(k-1)$-plane
of $\RP^{n+1}$, it follows by induction that

\Lemma{}{All spheres belonging to an elementary $(k+1)$-cell, $1\leq k\leq
m-1$, of the net $F$, or $\hat{F}$, are contained in a $k$-dimensional plane
in $\RP^{n+1}$.}

\noindent
This lemma provides a method to reconstruct an elementary $m$-cell of
the net $F$, or $\hat{F}$, from one point and all spheres $S_j$ belonging
to those $k$-cells, $k\geq2$, that contain the initial point --- or, from
these $k$-cells themselves\Footnote{As mentioned above, a sphere $S_j$
belonging to corresponding edges of $F$ and $\hat{F}$ can generically be
reconstructed from the four points. However, if only points of one of the
nets, $F$ or $\hat{F}$, are given the spheres cannot be uniquely reconstructed
--- but, for example, choosing them to be planes in a Euclidean space (i.e.\
$\hat{F}\equiv e_{\infty}$) resolves this ambiguity.} (cf.\ Fig.~\KCells):
since at least $(m-1)$ spheres belonging to {\it any} elementary $(m-1)$-cell
are known, the corresponding $(m-2)$-planes in $\RP^{n+1}$ can be constructed.
Intersecting these $(m-2)$-planes yields some of the missing spheres ---
which can now be used to construct all $(m-3)$-planes corresponding to the
$(m-2)$-cells. This construction can be continued until all spheres are
constructed as the intersection of the lines corresponding to 2-cells.
The missing points of the $m$-cell can then be obtained by reflecting known
points at the constructed spheres.

This construction scheme provides a method to solve  the following (cf.\CDS)

\Theorem{~(Cauchy problem)}{Knowing the $k$-dimensional subnets of an
$m$-dimensional discrete curvature line net $F$, $k\geq2$, that pass
through a point in all $({}^m_{\kern0.1em k})$ different directions,
or $k$-dimensional subnets of discrete curvature line nets $F$ and $\hat{F}$,
that pairwise envelope discrete Ribaucour $k$-sphere congruences, the
net $F$, resp.\ both nets $F$ and $\hat{F}$, can be reconstructed.}

\noindent
This theorem can as well be proved directly on the level of the points of
an elementary $m$-cell --- in the case $m=3$ of lowest dimension, it can be
reduced to Miguel's theorem \Berger\ (cf.\DoSa,\KoSch):
\Figure{\hbox{\epsfbox{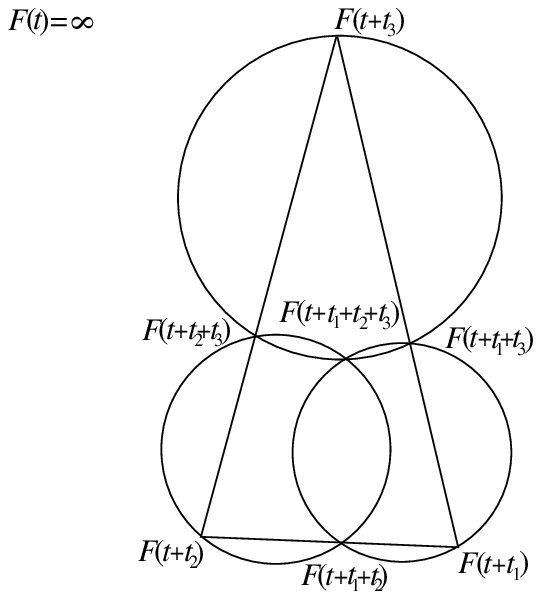}\vrule\epsfbox{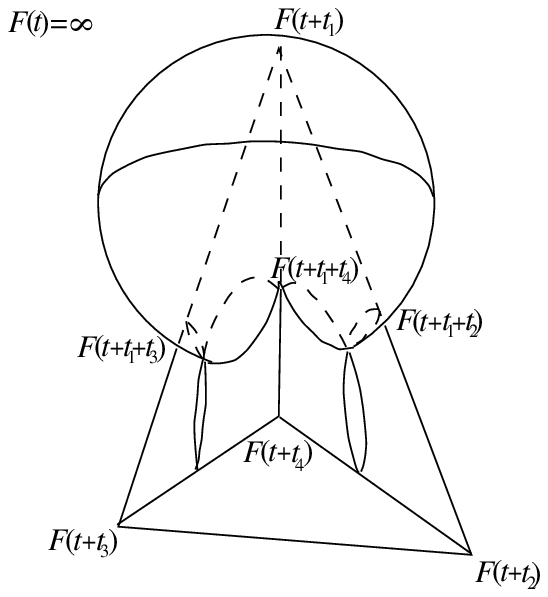}}}{Miguel's
  theorem}\Label\Miguel\noindent
since the vertices of an elementary 3-cell lie on a 2-sphere an inversion
that maps the common point $F(t)$ of all known 2-cells to infinity leaves
us with a 2-dimensional picture (Fig.~\Miguel, left). By Miguel's theorem
the three circles intersect in a point --- the image of the unknown point
of the elementary 3-cell.
A similar procedure works in higher dimensions: mapping the initial point
of an elementary $(k+1)$-cell to infinity, we are left with a $k$-simplex
in Euclidean $k$-space. Focussing on one of the $(k-1)$-spheres (corresponding
to the elementary $k$-cells) that contain the vertices of the simplex, the
problem can be reduced by one dimension (cf.\ Fig.~\Miguel, right).
Thus, the assumption follows by induction.

Another relation (besides relations through analogy) between the smooth
and the analogous discrete theory  of (Ribaucour pairs of) orthogonal
nets may be established through a

\Section{Permutability theorem}
for the Ribaucour transformation --- an orthogonal system $\hat{f}:
U\subset\R^m\to S^n$ (or, a submanifold with well defined curvature lines,
i.e.\ with flat normal bundle) is called  ``Ribaucour transformation'' of
a given orthogonal system $f:U\to S^n$ in the conformal $n$-sphere if they
form a Ribaucour pair.
``2-dimensional versions'' of this permutability theorem were known
classically (cf.eg.\Eisenhart):

\Theorem{}{Given two Ribaucour transforms $f_1,f_2:U\to S^3$ of a 2- or
3-dimensional orthogonal system $f:U\to S^3$ in the conformal 3-sphere,
there is a 1-parameter family of orthogonal systems $f_{12}=f_{21}:
U\to S^3$ that are Ribaucour transforms of both, $f_1$ as well as $f_2$.
Corresponding points of the four orthogonal systems are concircular.}

\noindent
Thus, a 2-dimensional discrete curvature line net $F:\Z^2\to S^3$ can be
obtained by repeatedly applying this permutability theorem\Footnote{Note,
that specializing this procedure to Darboux transforms of isothermic surfaces
one obtains discrete isothermic nets in a similar way (cf.\BoPi,\HJP).}
(cf.\KoSch): the images of one point in the original system will form a
discrete orthogonal (curvature line) net in $S^3$.
A ``3-dimensional version'' of this fact has recently been formulated \GaTsa: 

\Theorem{}{Given three Ribaucour transforms $f_1,f_2,f_3:U\to S^3$ of a
triply orthogonal system $f:U\to S^3$ as well as three systems
$f_{12},f_{23},f_{31}:U\to S^3$ such that each $f_{ij}$ is a Ribaucour
transform of $f_i$ and $f_j$ at the same time, there is a unique triply
orthogonal system $f_{123}:U\to S^3$ that is a Ribaucour transform of all
three, $f_{12}$, $f_{23}$ and $f_{31}$.}

\noindent
Clearly, this triply orthogonal system $f_{123}$ can be constructed as in
the solution of the Cauchy problem, earlier. In fact, starting with three
discrete curvature line nets, constructed according to the ``2-dimensional
version'' of the permutability theorem as described before, the solution of
the Cauchy problem provides a discrete triply orthogonal system in the
conformal 3-sphere.
Following this chain of thought, the Cauchy problem for discrete orthogonal
nets appears as a discrete analog of a general permutability theorem for
the Ribaucour transformation of smooth orthogonal systems --- or, the
permutability theorem appears as a smooth limit of the Cauchy problem:

\Conjecture{\Footnote{A rigorous proof seems to be a technicality.}}{Let
$1\leq i\leq k$ and let $f_t:U\subset\R^m\to S^n$ be smooth orthogonal
systems, corresponding to the vertices $t\in\{0,1\}^k$ of the
$({}^k_{\kern0.1em i})$ elementary $i$-cells containing $0\in\{0,1\}^k$,
such that the orthogonal systems corresponding to the endpoints of an edge
form a Ribaucour pair.
Then, all missing smooth orthogonal systems $f_t$, $t\in\{0,1\}^k$ in the
hypercube, can be --- uniquely if $i\geq2$ --- constructed such that any
edge corresponds to a Ribaucour pair.
Corresponding points of four orthogonal systems of an elementary 2-cell
are concircular.}

\noindent
Considering a smooth orthogonal system $f:U\subset\R^m\to\R^n$ as a discrete
orthogonal net with infinitesimal edge lengthes $f(t+\varepsilon t_i)-f(t)
\simeq0$ with fixed $\varepsilon\simeq0$:
the missing systems can be constructed by solving the Cauchy problem for
the $(m+k)$-dimensional\Footnote{Here, we may also allow $m+k>n$ without
harm.} net $F:\Z^{m+k}\to\R^n$ with the given smooth orthogonal systems as
part of $({}^{m+k}_{\kern0.8em i})$ initial data.
This might serve as a plausibility argument for the general permutability
theorem for the Ribaucour transformation of orthogonal systems.

\Acknowledgements{The major part of this paper was prepared in early 1997
when both authors stayed at the center of Geometry, Analysis, Numerics and
Graphics, GANG, at Amherst. We would like to thank the members of GANG and
of the mathematics department at the University of Massachusetts at Amherst
for their hospitality and their interest in our work.}

\References

\bye